\documentclass{article}

\usepackage{amssymb}
\usepackage{amsfonts}
\usepackage{pstricks}

\begin{document}

\title{Cluster networks and Bruhat--Tits buildings}

\author{S.V.Kozyrev\footnote{Steklov Mathematical Institute, Russian
Academy of Sciences}}

\maketitle

\begin{abstract}
Clustering procedure for the case where instead of a fixed metric one applies a family of metrics is considered. In this case instead of a classification tree one obtains a classification network (a directed acyclic graph with non directed cycles).

Relation to Bruhat--Tits buildings is discussed. Dimension of a general cluster system is considered.
\end{abstract}

Keywords: clustering, cluster networks, Bruhat--Tits buildings.

\section{Introduction}

This paper discusses the relation of construction of multidimensional cluster systems described in \cite{clustering}, \cite{cluster} and geometry of Bruhat--Tits buildings. Application to data analysis is also discussed.

Clustering procedure gives a construction of a tree of clusters with a hierarchy (partial order) starting from a metric on a set of points.
This procedure is an important method in data analysis, in particular in applications to bioinformatics  (construction of taxonomy, or a tree of life).

In the present paper we discuss the following approach to clustering: assume we have instead of one metric a family of metrics depending on a set of parameters (this is a typical situation in applications). We will get a family of clusterings. The question is: can we describe this family by a single mathematical object?

We will use an analogy from $p$-adic geometry. In $p$-adic spaces we have natural hierarchies (partially ordered trees) of balls. A hierarchy of this kind can be considered as a clustering with respect to metric in $p$-adic space.  In multidimensional $p$-adic geometry we have a generalization of hierarchies, described by the affine Bruhat--Tits buildings. These buildings are related to families of balls with respect to several ultrametrics.

In this paper we discuss a relation between Bruhat--Tits buildings and geometry of $p$-adic cluster systems with respect to a family of metrics. We also discuss a generalization of the corresponding structure of simplicial complex to general systems of clusters. In particular, we discuss a notion of dimension for general cluster systems.

For other application of $p$-adic numbers to data analysis see \cite{Khrennikov1}, \cite{Khrennikov2}. 

The structure of the present paper is as follows.

In Section 2 we recall the definition of clustering procedure and discuss some examples of cluster networks for a family of metrics.

In Section 3 we recall the definition of the affine Bruhat--Tits buildings and discuss the relation of this building and networks of clusters in $p$-adic spaces.

In Section 4 we discuss a structure of simplicial complex for general cluster networks and a definition of dimension form general cluster systems.

\section{Examples of cluster networks}

Let us recall the standard definition of clustering. We will use for simplicity the nearest neighbor clustering (we could also consider more general clustering algorithms). For general discussion of clustering see \cite{Murtagh}.

\medskip

Let $(M,\rho)$ be an arbitrary metric space.

A sequence of points $a=x_0,x_1,\dots,\allowbreak x_{n-1},x_n=b$ in
$(M,\rho)$ is called an $\varepsilon$-chain connecting two points
~$a$ and~$b$ if $\rho(x_k,x_{k+1})\le\varepsilon$ for all $0\le k <n$.

If there exists an $\varepsilon$-chain connecting
$a$ and~$b$ then~$a$ and~$b$ are $\varepsilon$-connected.

The chain distance between $a$ and $b$ is defined as
$d(a,b)=\inf(\varepsilon$: $a$,~$b$ $\varepsilon$-connected$)$.

This distance has all properties of ultrametric excluding non--degeneracy (i.e. non coinciding points can have zero chain distance). Therefore chain distance defines an ultrametric on the set of equivalence classes of points in $M$, where $a$, $b$ are in the same equivalence class when $d(a,b)=0$.

In particular if the initial metric $\rho$ is an ultrametric, then the corresponding chain distance $d$ coincides with $\rho$.

\medskip

A cluster $C(i,R)$ in a metric space
$(M,\rho)$ is a ball with the center~$i$
and radius $R$ with respect to the chain distance, i.e. the set $\{j\in M\colon  d(i,j)\le R \}$.

\medskip

A clustering of a metric space~$M$ is a cluster set, satisfying:

i) every element in~$M$ belongs to some cluster;

ii) for any pair $a$, $b$ of elements in~$M$ there exists a minimal cluster $\sup(a,b)$ containing both elements;

iii) for arbitrary embedded clusters $A\subset B$  every increasing
sequence of embedded clusters $\{A_i\}$,
$A\subset\dots\subset A_i\subset A_{i+1}\subset\dots\subset B$ is finite.

\medskip

Clustering procedure generates a partially ordered tree of clusters (dendrogram) in the following way:

i) vertices are clusters;

ii) partial order is given by inclusion of clusters;

iii) edge connects two clusters nested without intermediaries.

\medskip

Multidimensional generalization of clustering \cite{clustering}, \cite{cluster} is introduced as a generalization of clustering for the case of several metrics.
Assume we have a family of metrics (which is a standard situation for applications in data analysis), say this family is parameterized by a set of real parameters.   In this case instead of one cluster tree we will obtain a family of cluster trees. Some of clusters for the different trees can coincide as sets. Identifying these clusters for the different trees we will obtain a network of clusters. We will discuss the relation of cluster networks of this type and Bruhat--Tits buildings.

Before the discussion of general definition let us consider  examples of multidimensional clustering.

\medskip

\noindent{\bf Example 1}. The case of a set of three points $A$, $B$, $C$ in a two--dimensional
real plane $\mathbb{R}^2$ with the standard Euclidean metric. Parameters defining the metric are
coordinates of the points in the plane.

\medskip

Cluster tree ${\cal A}_1$ (see Fig. 1).
The cluster set contains $A$, $B$, $C$, $AB$, $ABC$
(vertices of the cluster tree), edges join the vertices in
accordance with the growth of the clusters, i.e. the cluster tree contains the edges
$$
(A, AB),\quad (B, AB),\quad (AB, ABC),\quad (C, ABC).
$$
Here we denote by $ABC$ the cluster containing $A$, $B$ and $C$.

\psset{xunit=1cm}
\psset{yunit=1cm}

\begin{figure}\label{fig1}
\begin{pspicture}(0,4)
   \psline[linewidth=2pt](0,1.5)(1,2.5)(2,1.5)
    \psline[linewidth=2pt](1,2.5)(2,3.5)(4,1.5)
   \rput(0,1){\bf A}
   \rput(2,1){\bf B}
   \rput(4,1){\bf C}
   \rput(1,3){\bf AB}
   \rput(2,4){\bf ABC}
   \rput(6,2){\bf A}
   \rput(8,2){\bf B}
   \rput(7,4){\bf C}
   \psellipse[fillcolor=lightgray](7,2)(1.5,0.7)
\end{pspicture}
\rput(0,0){\bf Fig. 1: Cluster tree ${\cal A}_1$}
\end{figure}

\medskip

Cluster tree ${\cal B}_1$ (see Fig. 2).
Let us consider a deformation of the metric (motion of the points in the plane $\mathbb{R}^2$) with the replacing of the above cluster set by the set of vertices $A$, $B$, $C$, $BC$, $ABC$ with the corresponding edges
$$
(B, BC),\quad (C, BC),\quad (BC, ABC),\quad (A, ABC).
$$

\begin{figure}\label{fig2}
\begin{pspicture}(0,4)
   \psline[linewidth=2pt](0,1.5)(2,3.5)(4,1.5)
    \psline[linewidth=2pt](2,1.5)(3,2.5)
   \rput(0,1){\bf A}
   \rput(2,1){\bf B}
   \rput(4,1){\bf C}
   \rput(3,3){\bf BC}
   \rput(2,4){\bf ABC}
   \rput(6,2){\bf B}
   \rput(8,2){\bf C}
   \rput(7,4){\bf A}
   \psellipse[fillcolor=lightgray](7,2)(1.5,0.7)
\end{pspicture}
\rput(0,0){\bf Fig. 2: Cluster tree ${\cal B}_1$}
\end{figure}

\medskip

Cluster network ${\cal C}_1$ (see Fig. 3).
This network is a union of the trees of clusters ${\cal A}_1$ and ${\cal B}_1$
(where we identify the clusters which coincide as sets).

The vertex set of the network ${\cal C}_1$ contains the clusters
$$
A,\quad B,\quad C,\quad AB,\quad BC,\quad ABC
$$
the set of edges of ${\cal C}_1$ contains
$$
(A, AB), (B, AB), (AB, ABC), (C, ABC),
$$
$$
(B, BC),\quad (C, BC),\quad (BC, ABC),\quad (A, ABC).
$$

The partial order of vertices is given by inclusion of clusters.

Cycles in this graph describe the different histories of the growth of clusters (growth with respect to the different metrics).

\begin{figure}\label{fig3}
\begin{pspicture}(0,4)
   \psline[linewidth=2pt](0,1.5)(0.5,2.5)(2,1.5)
   \psline[linewidth=2pt](0.5,2.5)(2,3.5)(4,1.5)
   \psline[linecolor=red,linewidth=2pt,linestyle=dashed](2,1.5)(3.5,2.5)(4,1.5)
   \psline[linecolor=red,linewidth=2pt,linestyle=dashed](3.5,2.5)(2,3.5)(0,1.5)
   \rput(0,1){\bf A}
   \rput(2,1){\bf B}
   \rput(4,1){\bf C}
   \rput(0,2){\bf AB}
   \rput(2,4){\bf ABC}
   \rput(4,3){\bf BC}
   \rput(6,2){\bf A}
   \rput(8,2){\bf B}
   \rput(8,4){\bf C}
   \psellipse[fillcolor=lightgray](7,2)(1.5,0.7)
   \psellipse[linecolor=red,linestyle=dashed](8,3)(0.7,1.5)
\end{pspicture}
\rput(0,0){\bf Fig. 3: Cluster network ${\cal C}_1$}
\end{figure}

\medskip

\noindent{\bf Example 2}.\quad Let us consider a set of four points $A$, $B$, $C$, $D$ located
in the plane $\mathbb{R}^2$ at the vertices of some quadrangle. We will have the following trees of clusters.

\medskip

Cluster tree ${\cal A}_2$ (see Fig. 4). Clustering with respect to the plane metric gives the clusters
$$
A,\quad B,\quad C,\quad D,\quad AB,\quad CD,\quad ABCD.
$$
The set of edges has the form
$$
(A, AB), (B, AB), (C, CD), (D, CD), (AB, ABCD), (CD, ABCD).
$$

\begin{figure}\label{fig4}
\begin{pspicture}(0,5)
   \psline[linewidth=2pt](0,1.5)(1,2.5)(2,1.5)
   \psline[linewidth=2pt](4,1.5)(5,2.5)(6,1.5)
   \psline[linewidth=2pt](1,2.5)(3,4.5)(5,2.5)
   \rput(0,1){\bf A}
   \rput(2,1){\bf B}
   \rput(4,1){\bf C}
   \rput(6,1){\bf D}
   \rput(1,3){\bf AB}
   \rput(5,3){\bf CD}
   \rput(3,5){\bf ABCD}
   \rput(7,2){\bf A}
   \rput(9,2){\bf B}
   \rput(7,4){\bf C}
   \rput(9,4){\bf D}
   \psellipse[fillcolor=lightgray](8,2)(1.5,0.5)
   \psellipse[fillcolor=lightgray](8,4)(1.5,0.5)
\end{pspicture}
\rput(0,0){\bf Fig. 4: Cluster tree ${\cal A}_2$}
\end{figure}

\medskip

Cluster tree ${\cal B}_2$ (see Fig. 5).
Deformation of the mentioned quadrangle gives the cluster set
$$
A,\quad B,\quad C,\quad D,\quad AC,\quad BD,\quad ABCD
$$
with the edges
$$
(A, AC), (C, AC), (B, BD), (D, BD), (AC, ABCD), (BD, ABCD).
$$

\medskip

Cluster network ${\cal C}_2$ (see Fig. 6) is the union of the trees ${\cal A}_2$ and ${\cal B}_2$ of clusters.
This network contains the unions of the vertex sets and the edges sets in ${\cal A}_2$ and ${\cal B}_2$.

\begin{figure}\label{fig5}
\begin{pspicture}(0,5)
   \psline[linewidth=2pt](0,1.5)(2,2.5)(4,1.5)
   \psline[linewidth=2pt](2,1.5)(4,2.5)(6,1.5)
   \psline[linewidth=2pt](2,2.5)(3,4.5)(4,2.5)
   \rput(0,1){\bf A}
   \rput(2,1){\bf B}
   \rput(4,1){\bf C}
   \rput(6,1){\bf D}
   \rput(1.5,3){\bf AC}
   \rput(4.5,3){\bf BD}
   \rput(3,5){\bf ABCD}
   \rput(8,2){\bf A}
   \rput(10,2){\bf B}
   \rput(8,4){\bf C}
   \rput(10,4){\bf D}
   \psellipse[fillcolor=lightgray](8,3)(0.5,1.5)
   \psellipse[fillcolor=lightgray](10,3)(0.5,1.5)
\end{pspicture}
\rput(0,0){\bf Fig. 5: Cluster tree ${\cal B}_2$}
\end{figure}

\begin{figure}\label{fig6}
\begin{pspicture}(0,5)
   \psline[linewidth=2pt](0,1.5)(1,2.5)(2,1.5)
   \psline[linewidth=2pt](4,1.5)(5,2.5)(6,1.5)
   \psline[linewidth=2pt](1,2.5)(3,4.5)(5,2.5)
   \psline[linecolor=red,linewidth=2pt,linestyle=dashed](0,1.5)(2,2.5)(4,1.5)
   \psline[linecolor=red,linewidth=2pt,linestyle=dashed](2,1.5)(4,2.5)(6,1.5)
   \psline[linecolor=red,linewidth=2pt,linestyle=dashed](2,2.5)(3,4.5)(4,2.5)
   \rput(0,1){\bf A}
   \rput(2,1){\bf B}
   \rput(4,1){\bf C}
   \rput(6,1){\bf D}
   \rput(1,3){\bf AB}
   \rput(5,3){\bf CD}
   \rput(1.9,3){\bf AC}
   \rput(4.1,3){\bf BD}
   \rput(3,5){\bf ABCD}
   \rput(8,2){\bf A}
   \rput(10,2){\bf B}
   \rput(8,4){\bf C}
   \rput(10,4){\bf D}
   \psellipse[fillcolor=lightgray](9,2)(1.5,0.5)
   \psellipse[fillcolor=lightgray](9,4)(1.5,0.5)
   \psellipse[linecolor=red,linestyle=dashed](8,3)(0.5,1.5)
   \psellipse[linecolor=red,linestyle=dashed](10,3)(0.5,1.5)
\end{pspicture}
\rput(0,0){\bf Fig. 6: Cluster network ${\cal C}_2$}
\end{figure}

\section{Affine Bruhat--Tits buildings and cluster networks}

In the present section we will show that for the network of balls in $\mathbb{Q}_p^d$ there exists a natural structure of simplicial complex which is related to the affine Bruhat--Tits building. Discussion of buildings one can find in particular in \cite{Garrett}. For discussion of $p$-adic geometry (in particular lattices) see \cite{Weil}.

\medskip

\noindent{\bf Affine Bruhat--Tits building.}\quad
Vertices of the building are equivalence classes of lattices.
A lattice in $\mathbb{Q}_p^d$ is an open compact $\mathbb{Z}_p$-module in $\mathbb{Q}_p^d$. Any lattice can be put in the form
$$
\oplus_{i=1}^{d}\mathbb{Z}_p e_i,
$$
where $\{e_i\}$ is a basis in $\mathbb{Q}_p^d$.

Two lattices are equivalent if one is a scalar multiple of the other.

Two lattices $L_1$ and $L_2$ are adjacent (connected by an edge) if some representatives from equivalence classes $L_1$ and $L_2$ satisfy
$$
pL_1\subset L_2\subset L_1.
$$

$k-1$-Simplices are defined as equivalence classes of $k$  adjacent lattices, i.e. the chains
$$
pL_{k}\subset L_1 \subset L_2\subset \dots \subset L_{k}.
$$
Here $1\le k\le d$.

An apartment in the affine building is the subcomplex corresponding to a fixed basis $\{e_i\}$ in $\mathbb{Q}_p^d$ which contains the equivalence classes of lattices $\oplus_{i=1}^{d}\mathbb{Z}_p p^{a_i}e_i$,  $a_i\in\mathbb{Z}$.

\medskip

\noindent{\bf Multidimensional $p$-adic metric.}\quad
Let us consider a metric $s_{q_1,\dots,q_d}(x,y)$ in $\mathbb{Q}_p^d$ defined by the norm $N_{q_1,\dots,q_d}(z)$
\begin{equation}\label{dq0}
s_{q_1,\dots,q_d}(x,y)=N_{q_1,\dots,q_d}(x-y),
\end{equation}
\begin{equation}\label{dq}
N_{q_1,\dots,q_d}(z)={\rm max}_{i=1,\dots,d}(q_i|z_i|_p),\qquad q_i \ne 0.
\end{equation}
Dilations $p^k\mathbb{Z}_p^d$, $k\in\mathbb{Z}$ are balls
with respect to all such norms $N_{q_1,\dots,q_d}$ if $p^{-1}<q_i\le 1$.

Here we use the following definition of norm for $\mathbb{Q}_p^d$: a norm is a function $N(\cdot)$ on $\mathbb{Q}_p^d$ taking values in $[0,\infty)$ and satisfying the conditions:

i) Nondegeneracy: $N(x)=0\Leftrightarrow x=0$;

ii) Linearity: $N(ax)=|a|_p N(x)$, $x\in\mathbb{Q}_p^d$, $a\in\mathbb{Q}_p$;

iii) Strong triangle inequality: $N(x+y)\le \max \left[N(x),N(y)\right]$.

A general norm ($A$-rotation of $N_{q_1,\dots,q_d}$) is defined as
\begin{equation}\label{sAq}
N^{A}_{q_1,\dots,q_d}(z)=N_{q_1,\dots,q_d}(Az),
\end{equation}
where $A$ is a matrix from ${\rm Gl}_d(\mathbb{Q}_p)$. A metric $s^{A}_{q_1,\dots,q_d}$ is defined by the norm $N^{A}_{q_1,\dots,q_d}$ as above (\ref{dq0}).


In particular for a norm $N_{q_1,\dots,q_d}$ of the form (\ref{dq}) with
\begin{equation}\label{condition}
p^{-1}<q_1<\dots<q_d\le 1
\end{equation}
the set of intermediary $N_{q_1,\dots,q_d}$--balls between $p\mathbb{Z}_p^d$ and $\mathbb{Z}_p^d$ contains the balls
$$
B_j=\mathbb{Z}_p\times\dots\times\mathbb{Z}_p\times p\mathbb{Z}_p\times\dots\times p \mathbb{Z}_p
$$
with $j$ components $\mathbb{Z}_p$ and $d-j$ components $p\mathbb{Z}_p$, $j=0,\dots,d$.

\medskip

\noindent{\bf Simplicial complex of balls.}\quad
Let us define a structure of simplicial complex on the network ${\cal C}$ of balls with respect to the defined above family of metrics $s^{A}_{q_1,\dots,q_d}$, $A\in{\rm Gl}_d(\mathbb{Q}_p)$.

Let $s$ be a metric from the described family and $I$ a $s$-ball containing zero (a $s$-ball is a ball with respect to $s$, zero is a vector in $\mathbb{Q}_p^d$ with zero coordinates). Then the dilation $pI$ is also a $s$-ball. The (containing zero) $s$-balls $I$ and $J$ are adjacent if $pI\subset J\subset I$.  $k-1$-Simplices are defined as families of $k$ adjacent $s$-balls
$$
pI_{k}\subset I_1 \subset I_2\subset \dots \subset I_{k}.
$$

Let us consider the maximal sequence of nested intermediary $s$-balls between $pI$ and $I$. If the parameters $q_i$ of the norm are generic (any two parameters can not be made equal by multiplication by degrees of $p$, for example when the parameters satisfy (\ref{condition})) then the above sequence contains $d+1$ balls and defines a $d-1$-simplex.

General simplices in the simplicial complex of balls with respect to a metric $s^{A}_{q_1,\dots,q_d}$ are defined as translations of simplices described above (translations as families of sets in $\mathbb{Q}_p^d$).

The simplicial complex ${\cal C}$ of balls with respect to the family $\{s^{A}_{q_1,\dots,q_d}\}$ of metrics is defined as a union of complexes of balls for different metrics  $s^{A}_{q_1,\dots,q_d}$. Here we identify $s$-ball and $s'$-ball which coincide as sets (and identify $s$ and $s'$-simplices which coincide as sets of balls).

\medskip

\noindent{\bf Relation between norms in $\mathbb{Q}_p^d$ and simplices in the affine building.}\quad Any ball with respect to a norm in $\mathbb{Q}_p^d$ which contains zero is a lattice. This follows from the strong triangle inequality.

Let us consider in the defined above simplicial complex ${\cal C}$ of balls the subcomplex ${\cal C}_0$ of balls which contain zero. For any ball $I$ containing zero a dilation $p^{k}I$, $k\in\mathbb{Z}$ is also a ball (with respect to the same norm). The same holds for simplices. Therefore the factor ${\cal C}_0/\Gamma$ by the group of dilations by $p^k$, $k\in\mathbb{Z}$ is a simplicial complex.

There exists a natural simplicial map from the simplicial complex ${\cal C}_0/\Gamma$ to the affine Bruhat--Tits building which put in correspondence to a ball the corresponding lattice.

The defined map is an embedding of the complex ${\cal C}_0/\Gamma$ into the affine building. Let us show that this map is surjective (i.e. is an isomorphism of simplicial complexes).

We say that two norms are equivalent if they generate the same family of balls. Let us show that to any maximal simplex in the affine Bruhat--Tits building one can put in correspondence an equivalence class of norms in $\mathbb{Q}_p^d$.

\medskip

Let $L$ be a lattice in $\mathbb{Q}_p^d$ and
\begin{equation}\label{lattices}
pL=L_0 \subset L_1\subset\dots\subset L_d=L
\end{equation}
be a maximal sequence of (different) embedded lattices (equivalently, a maximal simplex in the affine building).

For any pair $L_{j}\supset L_{j-1}$, $j=1,\dots,d$ of consecutive lattices in the above sequence let us choose an element $f_j\in L_{j}$, $f_{j}\notin L_{j-1}$. This gives a set $\{f_1,\dots,f_d\}$ of vectors in $\mathbb{Q}_p^d$. One has the following lemma.

\medskip

\noindent{\bf Lemma 1}\quad {\sl
1) The defined above set $\{f_1,\dots,f_d\}$ is a basis in $\mathbb{Q}_p^d$;

2) The lattices $L_j$ have the form
\begin{equation}\label{simplex}
L_j=\oplus_{i=1}^{j}\mathbb{Z}_pf_i\oplus \oplus_{i=j+1}^{d}p\mathbb{Z}_pf_i.
\end{equation}
}
\medskip

Let us introduce a norm in $\mathbb{Q}_p^d$ as follows. Let us put in correspondence to lattices $L_j$ from (\ref{lattices}) some positive numbers $q_j$,  $p^{-1}<q_1<\dots<q_d\le 1$.

Let us define a function $N(x)$ on $\mathbb{Q}_p^d$ in the following way. For $x\in L_j\backslash L_{j-1}$, $j=1,\dots,d$ we define $N(x)=q_{j}$. For $x=0$ we put $N(x)=0$. We define $N(\cdot)$ in all $\mathbb{Q}_p^d$ using the condition $N(p^{k}x)=p^{-k}N(x)$, $k\in\mathbb{Z}$.

\medskip

\noindent{\bf Lemma 2}\quad {\sl The function $N(\cdot)$ defined as above will be a norm in $\mathbb{Q}_p^d$ satisfying the strong triangle inequality. The sequence (\ref{lattices}) of lattices will be a maximal sequence of balls with respect to $N(\cdot)$ which lie between the balls $pL$ and $L$.}

\medskip

The introduced norm $N$ belongs to the family (\ref{sAq}) with the parameters satisfying (\ref{condition}). In particular, the matrix $A$ can be chosen as the matrix which maps the basis $\{f_1,\dots,f_d\}$ from the above lemma to the coordinate basis in $\mathbb{Q}_p^d$ and $N(\cdot)=N^{A}_{q_1,\dots,q_d}(\cdot)$.  Any two norms defined in this way will be equivalent (will generate the same set of balls).

We have constructed a norm of the form (\ref{sAq}), (\ref{condition}) starting from a maximal simplex in the affine building. Analogously, let us consider for a norm $N^{A}_{q_1,\dots,q_d}(\cdot)$ defined by (\ref{sAq}), (\ref{condition}) the set of lattices (\ref{simplex}) where the basis $\{f_1,\dots,f_d\}$ is defined by the matrix $A$ as above.
This set defines a simplex in the affine building.

We have shown that there exists a one to one correspondence between the equivalence classes of norms of the form (\ref{sAq}), (\ref{condition}) and maximal simplices in the affine Bruhat--Tits building.

In the above construction it is important to consider norms with generic parameters.
Let us consider a general norm $N$ of the form (\ref{dq}), (\ref{sAq}) and take a $N$-ball $L$.
It is possible (say if some parameters $q_i$ in (\ref{dq}) are equal) that a set of intermediary balls between $L$ and $pL$ contains less than $d+1$ balls and therefore can not define a $d-1$--simplex in the affine building.

\section{General networks of clusters}

In the present section we discuss a generalization to general cluster systems of the described in the previous sections cluster network ${\cal C}$ related to a family of metrics. We introduce a structure of simplicial complex on this network and a definition of dimension for cluster systems.

Let $X$ be a locally compact ultrametric space with some finite family of ultrametrics ${\bf s}$ defined on $X$. Moreover, let, for any pair of metrics $s,r\in{\bf s}$, any $s$-ball be a finite union of $r$-balls.

The family ${\bf s}$ of ultrametrics on $X$ is compatible, if for any two balls,  $s$-ball $I$ and $r$-ball $J$, $s,r\in{\bf s}$, the intersection $I\bigcap  J$ is a ball with respect to some ultrametric $t\in{\bf s}$.

We put in correspondence to a metric $s\in {\bf s}$ the corresponding tree ${\cal T}(X,s)$ of $s$-balls in $X$. Vertices of this tree are $s$-balls, two vertices are connected by edge if the corresponding balls are nested without intermediaries.

The graph ${\cal C}(X,{\bf s})$ (the network of clusters in $X$ with respect to the family ${\bf s}$ of metrics) is a union of trees ${\cal T}(X,s)$ of $s$-balls, $s\in{\bf s}$.

The set of vertices of ${\cal C}(X,{\bf s})$ is the union of the sets of $s$-balls, $s\in{\bf s}$, edges connect $s$-balls (with the same $s$) nested without intermediaries. The partial order in ${\cal C}(X,{\bf s})$ is defined by inclusion of subsets in $X$. If some $s$-ball coincides with some $r$-ball as a set, they define the same vertex in ${\cal C}(X,{\bf s})$.

Let us define simplices in the network ${\cal C}(X,{\bf s})$. Let ${\bf r}\subset {\bf s}$ be a subfamily of metrics in $X$.
Let us fix some ${\bf r}$-ball  $I$ (i.e. $I$ is a $s$-ball with respect to all $s\in{\bf r}$).
Let $J$ be a smallest ${\bf r}$-ball which is strictly greater than $I$.

We define a simplex as a subset (of cardinality at least two) of the set of intermediary $s$-balls lying between $I$ and $J$ (for metrics $s\in{\bf r}$). $k-1$-Simplex will contain $k$ $s$-balls, in particular, an edge in ${\cal T}(X,s)$ (a pair ball -- maximal subball) will be a one-simplex.

We have defined a structure of simplicial complex on the tree ${\cal T}(X,s)$. The union over metrics $s\in{\bf s}$ of these simplicial complexes defines a structure of simplicial complex on the network ${\cal C}(X,{\bf s})$ (where we, as usual, identify vertices which correspond to balls coinciding as sets and identify simplices coinciding as sets of vertices).

The idea of this definition of simplices is taken from the $p$-adic case discussed in the previous section. In this case simplices in the affine Bruhat--Tits building are given by sets of balls lying between balls $I$ and $pI$. If $I$ is a ball with respect to several metrics $s^{A}_{q_1,\dots,q_d}$ (say given by different rotations $A$ or different order of indices $q_j$) then $pI$ will also be a ball with respect to the same metrics. Therefore the above definition of simplicial complex ${\cal C}(X,{\bf s})$ generalizes the definition of the simplicial complex ${\cal C}$ of balls in $\mathbb{Q}_p^d$ of the previous section (instead of consideration of simplices as subsets of a sequence of nested balls between $I$ and the dilation $pI$ we consider cycles in the cluster network ${\cal C}(X,{\bf s})$, simplex is a subset of a path between the minimal and maximal vertices of a cycle).

\medskip

\noindent{\bf Dimension of clusters network.}\quad
Let us generalize the definition of dimension (the number of $p$-adic parameters) to general cluster networks.

Let ${\bf r}\subset {\bf s}$ be a subfamily of metrics in $X$. Let $I$, $J$ be a pair of ${\bf r}$-balls ($s$-balls for all $s\in {\bf r}$), where $I\subset J$ and $J$ is a minimal ${\bf r}$-ball which contains $I$.

Let us consider the maximal $s$-simplex, $s\in {\bf r}$ corresponding to $I$, $J$ (i.e. the maximal sequence of nested balls lying between $I$ and $J$).

The rank of this simplex we call the ${\bf r}$--dimension for the pair $I$, $J$\footnote{In principle this rank may differ for different $s\in{\bf r}$, in this case we take the maximal rank.}.

\medskip

In the $p$-adic case: the family ${\bf r}$ is given by the set of metrics $s^{A}_{q_1,\dots,q_d}$ with fixed matrix $A$ and all possible reorderings of parameters $p^{-1}< q_1<\dots<q_d\le 1$. The above dimension is equal to the number $d$ of $p$-adic parameters.

The introduced dimension is not equal to the VC (or Vapnik -- Cher\-vonenkis, or combinatorial) dimension \cite{thenature}.

\medskip

\noindent{\bf Applications to data analysis}.\quad
Clustering is a tool of data analysis with many applications, in particular to bioinformatics.
The set $X$ of data may be generated in a complex way, there may be some independent contributions in the data.

There should be some way to describe independencies in data at the level of networks of clusters.
Classification trees (trees of clusters) describe the diversity of data, the multidimensional generalization should describe the situation where we have independent sources of diversity.

Dimension of a network of clusters will describe the number of sources of diversity.

\medskip

\noindent{\bf {Applications to taxonomy for reticulate evolution}}.\quad
Let us discuss the application of the described above cluster networks in bioinformatics. One of the problems in bioinformatics is a construction of phylogenetic classification trees by comparison of genetic markers (some subsequences in genomes). These trees are constructed with the help of clustering procedure using some metric for genetic markers.

We are interested in the case when we have several genetic markers. In this case the metric is not uniquely defined.
This metric has the form of a sum of contributions from the different genetic markers
$$
d(X,Y)=\sum_{j=1}^N w_j d_j(X_j,Y_j),
$$
where $w_j\ge 0$ are weights, $X$ and $Y$ are genomes,
$X_j$ are $Y_j$ are genetic markers,
$d_j$ is the distance for the $j$-th genetic marker.

Different sets of weights generate the different classification trees. In particular, when only one weight
$d_j$ is non zero, this weight generates the classification tree for the corresponding genetic marker.
Union of these trees sometimes is called the forest of life \cite{Evo3}.

The observation is that for some genetic markers the corresponding classification trees are different because different genetic markers may have different evolution histories.
The evolution is {\it reticulate} --- some parts of a genome may have the different origin due to hybridization or horizontal gene transfer.
In this case instead of phylogenetic trees one can consider {\it phylogenetic networks}.
Mathematical methods of analysis of phylogenetic networks one can find in \cite{network1,network2}.

Discussed in the present paper classification networks might give a general framework for construction of phylogenetic networks. In this approach instead of reproduction of the detailed genetic history of populations (which in general is not possible) one could use classification networks for coarse grained description of the evolution of ensembles of genes.

\bigskip

\noindent{\bf Acknowledgments}\qquad
This paper was partially supported by the Program  ''Modern problems of theoretical mathematics'' of the Department of Mathematics of the Russian
Academy of Sciences.

\end{document}